\documentclass[12pt]{article}

\usepackage{graphicx}
\usepackage{amsmath}
\usepackage{amssymb}
\usepackage{latexsym}
\usepackage{subfigure}
\usepackage{crop}
\usepackage[section,subsection,subsubsection]{placeins}
\usepackage[colorlinks = true, pdfstartview = FitV, linkcolor = blue, citecolor = blue, urlcolor = blue]{hyperref}
\usepackage{bm}
\usepackage{enumerate}

\newcommand {\uu}  { {\bf u} }

\newcommand {\xx}  { {\bf x} }

\newcommand {\rr}  { {\bf r} }

\newcommand{\R}{{\rm I\!R}}

\newcommand {\veps} {\varepsilon}

\newcommand {\Ex} { {\mathbb E} }

\newcommand {\qq}  { {\bf q} }

\newcommand {\mm}  { {\bf m} }
\newcommand {\vv}  { {\bf v} }
\newcommand {\ww}  { {\bf w} }
\newcommand {\ff}  { {\bf f} }

\newcommand {\bb}  { {\bf b} }

\newcommand {\dd}  { {\bf d} }
\newcommand {\ee}  { {\bf e} }

\newcommand {\zero}  { {\bf 0} }

\newcommand {\beata} { { \boldsymbol \eta} }

\renewcommand{\vec}[1]{\ensuremath{\mathbf{#1}}}

\renewcommand{\div}{\nabla\cdot\,}
\newcommand{\grad}{\ensuremath {\vec \nabla}}

\newcommand{\defeq}{\mathrel{\mathop:}=}

\newtheorem{theorem}{Theorem}

\newcommand{\qed}{\nobreak \ifvmode \relax \else
      \ifdim\lastskip<1.5em \hskip-\lastskip
      \hskip1.5em plus0em minus0.5em \fi \nobreak
      \vrule height0.75em width0.5em depth0.25em\fi}

\begin{document}

\title{Algorithms that satisfy a stopping criterion, probably}

\author{Uri Ascher and Farbod Roosta-Khorasani
\thanks{Dept. of Computer Science, University of British Columbia, Vancouver, Canada
{\tt ascher/farbod@cs.ubc.ca} .
This work was supported in part by NSERC Discovery Grant 84306.}}

\maketitle

\begin{abstract}

Iterative numerical algorithms are typically equipped with a stopping criterion, where the
iteration process is terminated when some error or misfit measure is deemed to be
below a given tolerance.
This is a useful setting for comparing algorithm performance, among other purposes.

However, in practical applications a precise value for such a tolerance is rarely known; 
rather, only some possibly 
vague idea of the desired quality of the numerical approximation is at hand.
We discuss four case studies from different areas of numerical computation,
where uncertainty in the error tolerance value and in the stopping criterion is revealed in different ways.
This leads us to think of approaches to relax the notion of exactly satisfying a tolerance value.
 
We then concentrate on a {\em probabilistic} relaxation of the given tolerance.
This allows, for instance, derivation of proven bounds on the
sample size of certain Monte Carlo methods.
We describe an algorithm that becomes more efficient in a controlled way
as the uncertainty in the tolerance increases,
and demonstrate this in the context of some particular applications of inverse problems.           

\end{abstract}

%%%%%%%%%%%%%%%%%%%%%%%%%%%%%%%%%%%%%%%%%%%%%%%%%%%%%%%%%%%%%%%%%%%

%\begin{keywords} 
%error tolerance, mathematical software, iterative method,
%inverse problem, Monte Carlo method, trace estimation,
%large scale simulation, DC resistivity  
%\end{keywords}
%
%\begin{AMS}
%65C20, 65C05, 65M32, 65L05, 65F10
%\end{AMS}
%
%\pagestyle{myheadings}
%\thispagestyle{plain}
%\markboth{Satisfying a stopping criterion, probably}{Satisfying a stopping criterion, probably}

%%%%%%%%%%%%%%%%%%%%%%%%%%%%%%%%%%%%%%%%%%%%%%%%%%%%%%%%%%%%%%%%%%%%%%%%%%%%%

\section{Introduction}
\label{sec:int}

A typical iterative algorithm in numerical analysis and scientific computing requires a stopping criterion.
Such an algorithm involves a sequence of generated iterates or steps, an error tolerance, 
and a method to compute (or estimate) some quantity related to the error.
%Thus, one envisions a sequence of iterates generated by the algorithm
%together with an estimation method to determine some nonnegative quantity related to
%the error in a given iteration $k$, as well as a positive error tolerance.
If this error quantity is below the tolerance then the iterative procedure
is stopped and success is declared.

The actual manner in which the error in an iterate is estimated can vary
all the way from being rather complex
to being as simple as the normed difference between two consecutive iterates. Further,
the ``tolerance'' may actually be a set of values involving combinations of absolute and relative
error tolerances. There are several fine points to this, often application-dependent,
that are typically incorporated in mathematical software packages 
(see for instance {\sc Matlab}'s various packages for solving ordinary differential equation (ODE) 
or optimization problems). 
That makes some authors of introductory texts devote significant attention to the issue, 
while others attempt to ignore
it as much as possible (cf.~\cite{dabj,heath,agbook}). Let us choose here the middle way of      
considering a stopping criterion in a general form
\begin{eqnarray}
{\rm error\_estimate} (k) \leq \rho ,
\label{1}
\end{eqnarray}
where $k$ is the iteration or step counter, and $\rho > 0$ is the tolerance, assumed given.

But now we ask, {\em is $\rho$ really given?!} Related to this, we can also ask,
{\em to what extent is the stoppping criterion adequate?}

\begin{itemize}
\item
The numerical analyst would certainly {\em like} $\rho$ to be given.
That is because their job is to invent new algorithms, prove various
assertions regarding convergence, stability, efficiency, and so on,
and compare the new algorithm to other known ones for a similar task.
For the latter aspect, a rigid deterministic tolerance for a trustworthy
error estimate is indispensable.

Indeed, in research areas such as image processing where criteria
of the form~\eqref{1} do not seem to capture certain essential
features and the ``eye norm'' rules, 
a good comparison between competing algorithms can be far more delicate.
Moreover, accurate comparisons of algorithms that require stochastic input
can be tricky in terms of reproducing claimed experimental results. 
\item
On the other hand, a practitioner who is the customer of numerical algorithms,
applying them in the context of some
complicated practical application that needs to be solved,
will more often than not find it
very hard to justify a particular choice of a precise value for $\rho$ in~\eqref{1}.  
\end{itemize}

Our first task in what follows is to convince the reader that often
in practice there is a significant uncertainty in the actual
selection of a meaningful value for the error tolerance $\rho$, 
a value that must be satisfied.
Furthermore, numerical analysts are also subconsciously aware of this
fact of life, even though in most numerical analysis papers such a
value is simply given, if at all, in the numerical examples section.
%can be quick to not consider $\rho$ as 
%something less than 
%a ``holy constant'':
%we adapt to weaker conditions in different ways, depending on the situation
%and the advantage to be gained in relaxing the notion of an error tolerance.
%Three 
Four typical yet different classes of problems and methods are considered
in Section~\ref{sec:examples}.

Once we are all convinced that there is usually a considerable uncertainty
in the value of $\rho$ (hence, we only know it ``probably''), 
the next question is what to do with this notion.
The answer varies, depending on the particular application and the situation
at hand. In some cases, such as that of Section~\ref{sec2.1}, the effective
advice is to be more cautious, as mishaps can happen.
In others, such as that of Section~\ref{sec2.2}, we are simply led to
acknowledge that the value of $\rho$ may come from thin air (though one then concentrates
on other aspects). But there are yet other classes of
applications and algorithms, such as in Sections~\ref{sec2.3} and~\ref{sec2.4}, 
for which it makes sense to attempt to quantify
the uncertainty in the error tolerance $\rho$ using a probabilistic framework.
We are not proposing in this article to propagate an entire probability distribution for $\rho$:
that would be excessive in most situations. But we do show, by studying
an instance extended to a wide class of problems, that employing such a 
framework can be practical and profitable.

%\uri{The rest of this section is still to be revised!}

Following Section~\ref{sec2.4}
we therefore consider in Section~\ref{sec:prob} a particular manner of
relaxing the notion of a deterministic error tolerance, by allowing
an estimate such as~\eqref{1} to hold only within some given probability.
% limit.
%Thus we obtain in effect a pair of tolerances replacing $\rho$.
%
Relaxing the notion of an error tolerance in such a way allows
the development of theory towards an uncertainty quantification of
Monte Carlo methods (e.g., \cite{achlioptas,avto,yori,ipwe,hoip}). 
We concentreate on
%use as a case study 
the problem of estimating the trace of symmetric positive semi-definite (SPSD)
matrices that are given implicitly, i.e., as a routine to compute their product
with any vector of suitable size. We then use this to satisfy in a probabilistic sense
an error tolerance for approximating a data misfit function
for cases where many data sets are given.
Some details of results of this type that are relevant in the present context
and have been proved in~\cite{roas1,roszas} 
%and (Roosta-Szekely-Ascher) 
are provided.

In Section~\ref{sec:application} we 
%quickly 
apply the results of Section~\ref{sec:prob}
%for implicit matrix trace estimation 
to form a randomized
algorithm for approximately solving inverse problems involving 
partial differential equation (PDE) systems.
We concentrate on cases where many data sets are provided 
(which is typical in several important applications) that
correspondingly require the solution of many PDEs.   
Conclusions and some additional general comments are offered in Section~\ref{sec:conclusions}.

While this article concentrates on bringing to the fore a novel point of view,
we also note that the observations in Section~\ref{sec2.2.1}, 
and to a lesser extent also Section~\ref{sec2.3.1}, 
are novel. The results displayed in Figures~\ref{fig0} and~\ref{fig02} are new.

%%%%%%%%%%%%%%%%%%%%%%%%%%%%%%%%%%%%%%%%%%%%%%%%%%%%%%%%%%%%%%%%%%%%%%%%%%%

\section{Case studies}
\label{sec:examples}

In this section we consider four classes of problems and associated algorithms,
in an attempt to highlight the use of different tests of the form~\eqref{1}
and in particular the implied level of uncertainty in the choice of $\rho$.

\subsection{Stopping criterion in initial value ODE solvers}
\label{sec2.1}

Using a case study, we show in this section that numerical analysts, too,
can be quick to not consider $\rho$ as 
%something less than 
a ``holy constant'':
we adapt to weaker conditions in different ways, depending on the situation
and the advantage to be gained in relaxing the notion of an error tolerance.

Let us consider an initial value ODE system in ``time'' $t$, written as
%first a boundary value ODE system written as
\begin{subequations}
\begin{eqnarray}
 \frac{d\uu}{dt} &=& \ff(t,\uu), \quad 0 \leq t \leq b, \label{1.1a} \\
\uu(0) &=& \vv_0 , \label{1.1b}
% \gb( \uu(0), \uu(b)) &=& \zero , \label{1.1b}
 \end{eqnarray}
 \label{1.1}
\end{subequations}
with $\vv_0$ a given initial value vector.
A typical adaptive algorithm proceeds to generate pairs
$(t_i,\vv_i), \ i = 0,1,2, \ldots, N$, in $N$ consecutive steps,
thus forming a mesh $\pi$ such that 
\[ \pi: 0 = t_0 < t_1 < \cdots < t_{N-1} < t_N = b ,\]
and $\vv_i \approx \uu(t_i), \ i = 1, \ldots , N$.

Denoting the numerical solution on the mesh $\pi$ by $\vv^\pi$, and
the restriction of the exact ODE solution to this mesh by $\uu^\pi$,
there are two general approaches for controling the error in such
an approximation.
\begin{itemize}
\item
%A typical procedure in a package for solving such a  problem numerically
%(see, e.g., \cite{amr} and references therein) 
%proceeds as follows:
Given a tolerance value $\rho$, keep estimating the {\em global error} and 
refining the mesh (i.e., the gamut of step sizes)
until roughly
\begin{eqnarray}
 \| \vv^\pi - \uu^\pi \|_\infty \leq \rho .
 \label{1.2}
 \end{eqnarray} 
Details of such methods can be found, for instance, in~\cite{hanowa,Higham91,cape,apbook}.
% where $\uu_\pi$ is the unknown exact solution restricted to the mesh or the approximate solution.
%(Matlab's, PASVA, ACR, Russell-Huang, ).

In \eqref{1.2} we could replace the absolute tolerance 
by a combination of absolute and relative tolerances, 
perhaps even different ones for different ODE equations.
But that aspect is not what we concentrate on in this article.

%Next consider an initial value ODE system, namely, a similar ODE as in \eqref{1.1a}
%but subject to an initial condition
%\begin{eqnarray*}
% \uu(0) &=& \vv_0 ,
%\end{eqnarray*}
%which is a special case of~\eqref{1.1b}.
%We could apply a similar procedure for stopping the iteration, namely, attempt to approximately satisfy
%\eqref{1.2} for a given $\rho$ \cite{hanowa,Higham91,cape}. 

\item
However, most general-purpose ODE codes estimate a {\em local error} measure for~\eqref{1}
instead, and refine  the step size locally. Such a procedure advances 
one step at a time, and estimates the next step size using local information
related to the local truncation error, or simply the difference between two
approximate solutions for the next time level, one of which presumed to be significantly more
accurate than the other.\footnote{
Recall that the local truncation error at some time $t = t_i$ is the amount by
which the exact solution $\uu^{\pi}$ fails to satisfy the scheme that defines $\vv^{\pi}$
at this point.
Furthermore, if at $t_i$, using the known $\vv_i$ and a guess for $t_{i+1}$,
we apply one step of two different Runge-Kutta methods of orders $4$ and $5$, say, then
the difference of the two results at $t_{i+1}$ gives an estimate for the error in
the lower order method over this mesh subinterval.} 
For details see~\cite{hanowa,hawa,apbook} and many references therein.
In particular, the popular {\sc Matlab} codes {\tt ode45} and {\tt ode23s}
use such a local error control.
\end{itemize}

The reason for employing local error control is that this allows for developing 
a much cheaper and yet
more sensitive adaptive procedure, an advantage that cannot be had, for instance, for 
general boundary value ODE problems; see, e.g.,~\cite{amr}.

But {\em does this always produce sensible results?!}
The answer to this question is negative. 
%\begin{itemize}
%\item
A simple example to the
contrary is the problem 
\[ \frac{du}{dt} =  100 (u - \sin t) + \cos t, \quad u(0) = 0, ~b = 1 . \]
Local truncation (or discretization) 
errors for this unstable initial value ODE propagate like $\exp (100 t)$,
a fact that is not reflected in the local behaviour of the exact solution $u(t) = \sin t$
on which the local error control is based.
Thus, we may have a large error $\| \vv^\pi - \uu^\pi \|_\infty$ even if
the local error estimate is bounded by $\rho$ for a small value of $\rho$.
%\item

%%%%%%%%%%%%%%%%%%%%%%%%%%%%%%%%%%%%%%%%%%%%%%%%%%%%%%%%%%%%

\subsubsection{Local error control can be dangerous even for a stable ODE system}
\label{sec2.1.1}

Still one can ask, {\em are we safe with local error control in case that we know
that our ODE problem is stable?}
Here, by ``safe'' we mean that the global error will not be much larger
than the local truncation error in scaled form.
The answer to this more subtle question turns out to be negative as well.
The essential point is that the global error consists of an accumulation
of contributions of local errors from previous time steps.
If the ODE problem is asymptotically stable (typically, because it describes a
damped motion) then local error contributions die away as time increases, 
often exponentially fast, so at some fixed time only the
most recent local error contributions dominate in the sum of contributions that forms
the global error. However, if the initial value ODE problem is merely marginally
stable (which is the case for Hamiltonian systems)
then local error contributions propagate undamped, and their
accumulation over many time steps can therefore be significantly larger 
than just one or a few such errors.\footnote{The local error control basically seeks
to equalize the magnitude of such local errors at different time steps.}  

%A less obvious (and less contrived) instance arises when 
For a simple concrete example, consider applying {\tt ode45}
with default tolerances to find the linear oscillator with a slowly varying frequency
that satisfies the following initial value ODE for $p(t)$:
\begin{eqnarray*}
 \frac{dq}{dt} &=& \lambda^2 p, \quad q(0) = 1, \\
 \frac{dp}{dt} &=& -(1+t)^2 q , \quad p(0) = 0 .
\end{eqnarray*}
Here $\lambda > 0$ is a given parameter.
Thus, $\uu = (q,p)^T$ in the notation of \eqref{1.1}.
This is a Hamiltonian system, with the Hamiltonian function given by
\[ H(q,p,t) = \frac 12 \left[ ((1+t)q)^2 + (\lambda p)^2 \right] .\]
%Thus, $\uu = (q,p)^T$ in the notation of \eqref{1.1}.
%, where $\frac{dp}{dt} = -\lambda^2 q$.
Now, since the ODE is not autonomous, the Hamiltonian is not constant in time.
However, the {\em adiabatic invariant}
\[ J(q,p,t) = H(q,p,t)/(1+t) \] 
(see, e.g.,~\cite{lere,asre}) is almost constant for large $\lambda$, satisfying
\[ [J(t) - J(0)]/J(0) = {\mathcal O} (\lambda^{-1})  \]
over the interval $[0,1]$.
%; see~\cite{lere,asre}.
This condition means in particular that for $\lambda \gg 1$ and the initial values
given above, 
$J(1) = J(0) + {\mathcal O} (\lambda^{-1}) \approx J(0)$. 
%is almost constant on this time interval.
%an exponentially large interval in $t$; see~\cite{asre} and references therein.

\begin{figure}[htb] 
\centerline{
{\includegraphics[width=.80\linewidth]{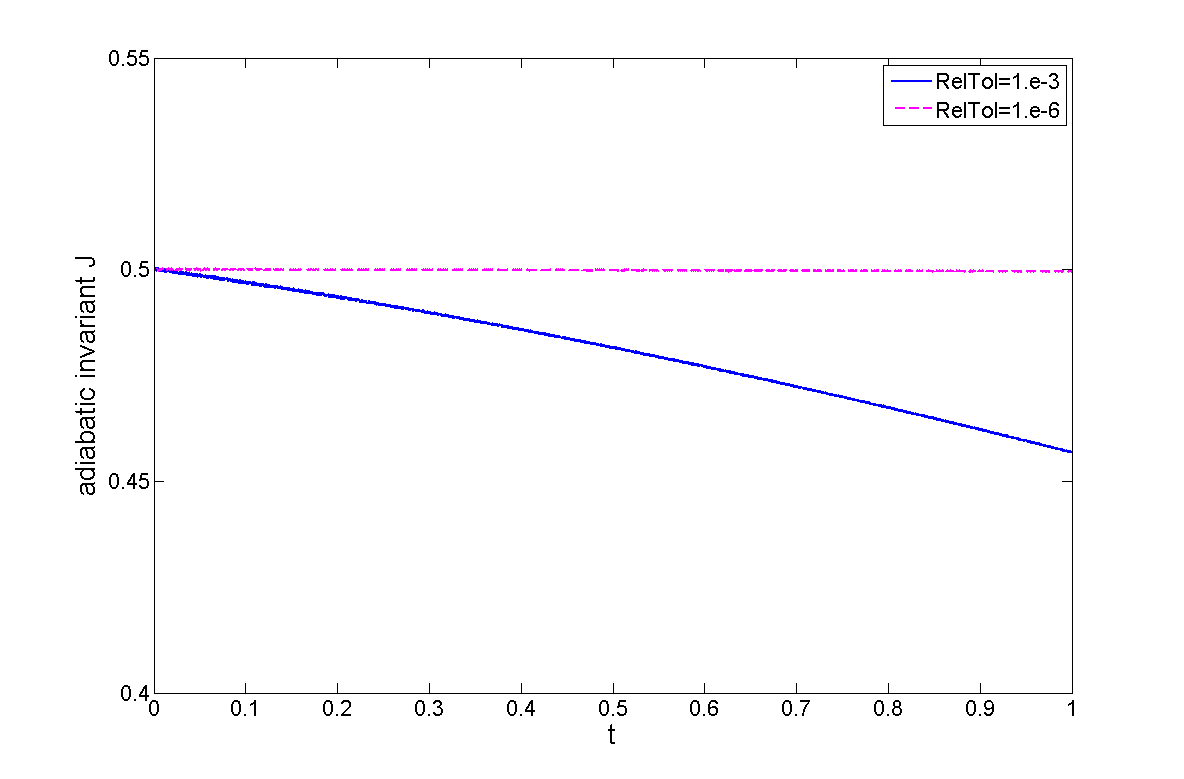}}
}
\caption{Adiabatic invariant approximations obtained using {\sc Matlab}'s package {\tt ode45}
with default tolerances (solid blue) and stricter tolerances (dashed magenta).
\label{fig0} }
\end{figure}
Figure~\ref{fig0} depicts two curves approximating the adiabatic invariant 
%overloading the notation $J(t)$, 
for $\lambda = 1000$. 
Displayed are the calculated curve using {\tt ode45} with default
tolerances (absolute=1.e-6, relative=1.e-3), as well as what is obtained upon
using {\tt ode45} with the stricter relative tolerance RelTol=1.e-6.
%(Note the hopefully unambiguous overloading of the notation $J(t)$ in the plot.)
From the figure it is clear that when using the looser
tolerance, the resulting approximation for $J(1)$ differs from $J(0)$ by 
far more than what $\lambda^{-1}=$1.e-3 and RelTol=1.e-3 
%or the relative tolerance 1.e-3  
would indicate, while the stricter tolerance gives a qualitatively correct result,
using the ``eye norm''.
Annoyingly, the qualitatively incorrect result does not look like ``noise'': 
while not being physical, it looks 
downright plausible, and hence could be misleading for an unsuspecting user.
Adding to the pain is the fact that this occurs for default tolerance values,
an option that a vast majority of 
%simple-minded 
users would automatically select.
$\blacklozenge$

A similar observation holds when trying to
approximate the phase portrait or other properties of an autonomous Hamiltonian ODE system
over a long time interval using {\tt ode45} with default tolerances: 
this may produce qualitatively wrong results.
See for instance Figures~16.12 and~16.13 in~\cite{agbook}:
the Fermi-Pasta-Ulam problem solved there is described in detail 
in Chapter~1 of~\cite{haluwa}.
What we have just shown here is that the phenomenon can arise also for a very
modest system of {\em two linear ODEs that do not satisfy any exact invariant}. 
%\end{itemize}   

We hasten to add that the documentation of {\tt ode45} (or other such codes) does not propose
to deliver anything like~\eqref{1.2}. Rather, the tolerance is just a sort of a knob 
that is turned to control
local error size. However, this does not explain the popularity of such codes despite
their limited offers of assurance in terms of qualitatively correct results.

Our key point in the present section is the following:
we propose that one reason for the popularity of ODE codes that use only
local error control is that in applications
one rarely knows a precise value for $\rho$ as used in~\eqref{1.2} anyway.
(Conversely, if such a global error tolerance value is known and is important
then codes employing a global error control, and not {\tt ode45}, should be used.)
Opting for local error control over global error control
 %Replacing the global error control by a local error control 
can therefore be seen as one specific way of adjusting 
%codes 
mathematical software in a deterministic sense to realistic uncertainties 
regarding the desired accuracy.

%%%%%%%%%%%%%%%%%%%%%%%%%%%%%%%%%%%%%%%%%%%%%%%%%%%%%%%%%%%%%%%%%%%%%%%%%%%

\subsection{Stopping criterion in iterative methods for linear systems}
\label{sec2.2}

In this case study, extending basic textbook material, we argue not only that
tolerance values used by numerical analysts are often determined solely for 
the purpose of the comparison
of methods (rather than arising from an actual application), but also that this
can have unexpected effects on such comparisons. 
In Section~\ref{sec2.2.1} we further make some novel, and we think intriguing,
observations on PDE discretizations, which arise in the present context although having
little to do with our main theme.

Consider the problem of finding $\uu$ satisfying
\begin{eqnarray}
 A \uu = \bb , 
\label{2.1}
\end{eqnarray}
where $A$ is a given $s \times s$ symmetric positive definite matrix such that
one can efficiently carry out matrix-vector products $A\vv$ for any suitable vector $\vv$,
but decomposing the matrix directly (and occasionally, even looking at its elements)
is too inefficient and as such is ``prohibited''.
We relate to such a matrix as being given {\em implicitly}.
The right hand side vector $\bb$ is given as well.

An iterative method for solving~\eqref{2.1} generates a sequence of iterates
$\uu_1, \uu_2, \ldots , \uu_k, \ldots$ for a given initial guess $\uu_0$.
Denote by $\rr_k = \bb - A\uu_k$ the residual in the $k$th iterate.
%It is well-known that 
The MINRES method, or its simpler version  Orthomin(2), 
%\cite{greenbaum} 
%conjugate gradient (CG) method [or orthomin? a cousin of minres;
%check; Greenbaum, my writeup]
can be applied to reduce the residual norm so that 
\begin{eqnarray}
\| \rr_k \|_2 \leq \rho \| \rr_0 \|_2
\label{2.2}
\end{eqnarray}
in a number of iterations $k$ that in general is at worst $\mathcal{O} \left(\sqrt{\kappa(A)}\right)$,
where $\kappa(A) = \| A \|_2 \|A^{-1}\|_2$ is the condition number of the matrix $A$.
Below in Table~\ref{table01} we refer to this method as MR.
The more popular conjugate gradient (CG) method generally performs comparably in practice.
We refer to~\cite{greenbaum} for the precise statements of convergence bounds and their proofs.

A well-known and simpler-looking family of {\em gradient descent} methods is given by 
%\begin{subequations}
\begin{eqnarray}
 \uu_{k+1} = \uu_k + \alpha_k \rr_k, 
\label{2.3}
\end{eqnarray}
where the scalar $\alpha_k > 0$ is the step size.
Such methods have recently come under intense scrutiny because of applications
in stochastic programming and sparse solution recovery.
Thus, it makes sense to evaluate and understand them in the simplest context of~\eqref{2.1},
even though it is commonly agreed that for the strict purpose of solving~\eqref{2.1} iteratively,
CG cannot be significantly beaten.  
Note that~\eqref{2.3} can be viewed as a forward Euler discretization of the artificial time ODE
\begin{eqnarray}
 \frac {d\uu}{dt} = -A\uu + \bb,\label{2.35}
\end{eqnarray}
with ``time'' step size $\alpha_k$. Next we consider two choices of this step size.

The steepest descent (SD) variant of~\eqref{2.3} is obtained by the
greedy (exact) line search for the function
\[ f(\uu) = \frac 12  \uu^TA\uu - \bb^T\uu , \]
which gives
\[ \alpha_k = \alpha_k^{SD} = \frac{\rr_k^T\rr_k}{\rr_k^TA\rr_k} \equiv \frac {(\rr_k, \rr_k)}{(\rr_k, A\rr_k)} 
\equiv \frac{\| \rr_k\|_2^2}{\| \rr_k\|_A^2} . \]
However, SD is very slow, requiring $k$ in \eqref{2.2} to
be proportional to $\kappa (A)$; see, e.g.,~\cite{akaike}.\footnote{
The precise statement of error bounds for CG and SD in terms of the error $\ee_k = \uu - \uu_k$
uses the $A$-norm, or ``energy norm'', and reads
\begin{eqnarray*}
\| \ee_k \|_A &\leq& 2 \left(\frac{\sqrt{\kappa(A)}-1}{\sqrt{\kappa(A)}+1}
\right)^k \|\ee_0\|_A, \quad {\rm for~CG}, \\
\| \ee_k \|_A &\leq& \left(\frac{\kappa(A)-1}{\kappa(A)+1}
\right)^k \|\ee_0\|_A, \quad {\rm for~SD}.
\end{eqnarray*}
See~\cite{greenbaum}.
}

A more enigmatic choice in~\eqref{2.3} is the lagged steepest descent
(LSD) step size
\[ \alpha_k = \alpha_k^{LSD} = \frac {(\rr_{k-1}, \rr_{k-1})}{(\rr_{k-1}, A\rr_{k-1})} . \]
It was first proposed in~\cite{babo} and practically used for instance in~\cite{befr08,dafl}.
To the best of our knowledge,
there is no known a priori bound on how many iterations as a function of $\kappa (A)$
are required to satisfy~\eqref{2.2} with this method~\cite{babo,RaydanThesis1991,fmmr,doas3}.

We next compare these four methods in a typical fashion for a 
prototypical PDE example, where we 
%this is all textbook material (e.g. ascher-greif).
consider the model Poisson problem
\[ -\Delta u =  1, \quad 0 < x,y < 1,\]
subject to homogeneous Dirichlet BC,
and discretized by the usual 5-point difference scheme
on a $\sqrt{s} \times \sqrt{s}$ uniform mesh.
%Set $s = N^2$ and 
Denote the reshaped vector of mesh unknowns by $\uu \in \R^{s}$. 
%The resulting condition number $\kappa(A)$
%of the matrix $A$ in~\eqref{2.1} is then proportional to $s$; see, e.g.,~\cite{agbook}.
%To be precise, 
The largest eigenvalue of the resulting matrix $A$ in~\eqref{2.1} is 
$\lambda_{\max} = 4h^{-2}(1 + \cos(\pi h))$,
and the smallest is $\lambda_{\min} = 4h^{-2}(1 - \cos(\pi h))$, 
where $h = 1 / (\sqrt{s}+1)$. Hence by Taylor expansion of $\cos(\pi h)$,
for $h \ll 1$ the condition number is essentially proportional to $s$:
\[ \kappa(A) = \frac{\lambda_{\max}} {\lambda_{\min}}  \approx \left(\frac 2{\pi} \right)^2 s . \]

In Table~\ref{table01} we list iteration counts required to satisfy~\eqref{2.2}
with $\rho = 10^{-7}$, starting with $\uu_0 = \zero$.

\begin{table}[!ht]
\begin{center}
\begin{tabular}{|rcccc|}
\hline
$s$ & MR & CG & SD & LSD  \\
\hline
$7^2$ & 9 & 9 & 196 & 45 \\
$15^2$ & 26 & 26 & 820 & 91 \\
$31^2$ & 54 & 55 & 3,337 & 261\\
$63^2$ & 107 & 109 & 13,427 & 632 \\
$127^2$ & 212 & 216 & 53,800 & 1,249 \\
\hline \hline
\end{tabular}
\end{center}
\caption{Iteration counts required to satisfy~\eqref{2.2}
for the Poisson problem with tolerance $\rho = 10^{-7}$ and  different mesh sizes $s$.
\label{table01}}
\end{table}

But now, returning to the topic of the present article, we ask,
{\em why insist on {$\rho = 10^{-7}$}?}
Indeed, the usual observation that one draws from the columns of values for MR, CG and SD
in a table such as Table~\ref{table01},
is that the first two grow like $\sqrt{\kappa(A)} \propto \sqrt{s}$  
while the latter grows like $\kappa(A) \propto s$.
The value of $\rho$, so long as it is not too large, does not matter at all!

\begin{figure}[htb] 
\centerline{
\mbox{\subfigure[Residuals]
{\includegraphics[width=.52\linewidth]{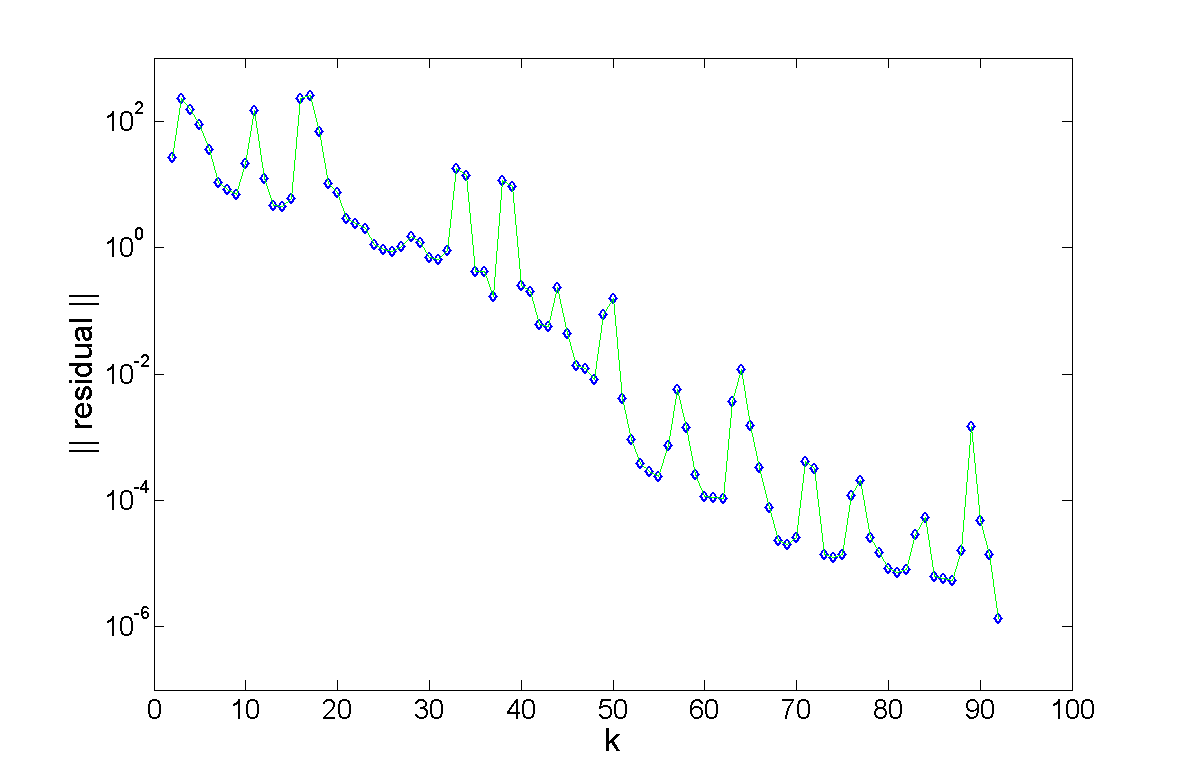}} \
\subfigure[Step sizes]
{\includegraphics[width=.52\linewidth]{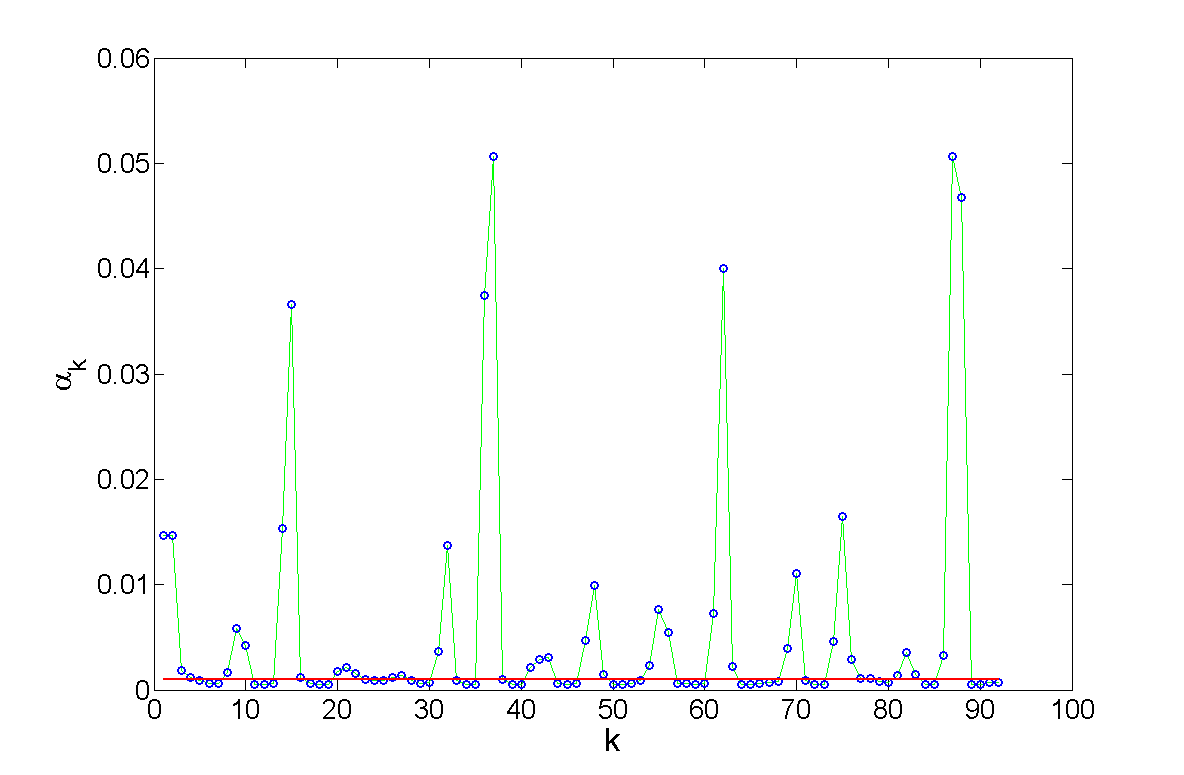}}
}}
\caption{Relative residuals and step sizes for solving the model Poisson problem
using LSD on a $15 \times 15$ mesh. The red line in (b) is the forward Euler stability limit.
\label{fig01} }
\end{figure}

And yet, this is not quite the case for the LSD iteration counts. These do not decrease in the
same orderly fashion as the others, even though they are far better 
(in the sense of being significantly smaller) than those for SD.
Indeed, this method is chaotic~\cite{doas3},
and the residual norm decreases rather non-monotonically, see Figure~\ref{fig01}(a).
Thus, the iteration counts in Table~\ref{table01} 
correspond to the iteration number $k=k^*$ where the rough-looking relative residual norm 
first records a value below the tolerance $\rho$. Unlike the other three methods, here the particular
value of the tolerance, picked just to be concrete, does play an unwanted role in the relative values,
 as a function of $s$, or $\kappa (A)$, of the listed iteration counts. $\blacklozenge$

\subsubsection{How large can the step size be? -- benefitting from a discretize-first approach}
\label{sec2.2.1}
 
While on the subject of LSD, let us also remark that gradient descent
for this example can be viewed as an application of the forward Euler method via~\eqref{2.35}
to find the steady state of the heat equation 
\begin{eqnarray}  
\frac{\partial u}{\partial t} = \Delta u + 1. \label{2.45}
\end{eqnarray}
(Homogeneous Dirichlet boundary conditions for all $t > 0$
are assumed already incorporated into the Poisson-Laplace operator $\Delta$.)

The forward Euler time-stepping method with step size $\alpha$, applied as a semi-discretization in time
to this PDE problem, reads
\begin{eqnarray}
 u_{k+1} = u_k + \alpha (\Delta u_k + 1), \label{2.46}
\end{eqnarray}
where $u_k = u_k(x,y)$.
However, this cannot be stable for any $\alpha > 0$. To see this, think of what happens
to a discontinuous initial condition function $u_0(x,y)$ when its spatial derivatives are repeatedly taken
as in~\eqref{2.46}. 

Next, let us discretize in space first, obtaining~\eqref{2.35}, and follow this with a forward
Euler discretization in time (i.e., apply gradient descent~\eqref{2.3}).
It is well-known that if the step size $\alpha_k = \alpha$ is constant for all $k$ 
then stability requires the restriction $\alpha \leq .25/s$.
In fact, the precise stability statement is $\alpha \leq 2/\lambda_{\max}$, where
$\lambda_{\max}$ is the largest eigenvalue of the matrix $A$ given above.
This restriction becomes more severe as $s$ grows, and in the limit $s \rightarrow \infty$
we obtain the unconditionally unstable method~\eqref{2.46}. 

So far, things are as expected.
But now, if $\alpha_k$ is allowed to vary with $k$ then no such stability requirement holds for all $k$!
Indeed, if $\alpha_k = 1/\lambda_{\min}$ then although this does enlarge the magnitude 
of high-frequency components in the residual (we are using multigrid jargon here; see~\cite{tos}),
it may decrease low frequency ones rapidly. The temporary increase in high frequency
component magnitudes can then be alleviated by next taking a few small step sizes which are
known to be particularly efficient in reducing those high frequency amplitudes.
Next we observe that as the uniform mesh size $s$ grows, $\lambda_{\min} \rightarrow 2\pi$,
which is independent of the discretization mesh width.
So this large step size is in effect constant as the mesh is refined, and yet
the method can be convergent.
This argument follows from a discretize-first approach, and cannot be seen from~\eqref{2.46}.

Returning to LSD as a method for automatically choosing the step sizes in the forward Euler
method, the results are recorded in Figure~\ref{fig01}(b). 
In fact, in this particular example we find upon decreasing the spatial mesh width  
%and keeping $\rho$ small
that $\max_k \alpha_k \approx .05$, independently of $s$.
Thus, what was merely possible in the previous paragraph is in fact automatically achieved
by an LSD method that provably converges to steady state. 

To the best of our knowledge, this result is novel, and moreover we have found
that it is unintuitive to many experts in numerical PDEs.
Hence we chose to incorporate it in this article, despite the fact that it
does not directly relate to its main theme.
%This effect does not seem to be as well-known in the community as perhaps it should be.

%%%%%%%%%%%%%%%%%%%%%%%%%%%%%%%%%%%%%%%%%%%%%%%%%%%%%%%%%%%%%%%%%%%%%%%%%%%

\subsection{Data fitting and inverse problems}
\label{sec2.3}

In the previous two case studies we have encountered situations where the intuitive use of
an error tolerance within a stopping criterion could differ widely (and wildly) from the notion that
is embodied in~\eqref{1} for the consumer of numerical analysts' products. 
We next consider a family of problems where the value of $\rho$ in a particular criterion~\eqref{1}
is more directly relevant.  

Suppose we are given observed data {$\dd \in {\R}^{l}$} and a {\em forward operator}
%$f(m)$,
{$f_i(m), \ i = 1, \ldots , l$,} 
which provides {predicted data} for each instance
of a distributed parameter function $m$. 
The (unknown) function $m$ is defined in some domain $\Omega$ in physical space and possibly time. 
We are particularly interested here in problems
where $f$ involves the solution $u$ in $\Omega$ of some linear PDE system,
sampled in some way at the points where the observed data are provided. 
Further, for a given mesh $\pi$ discretizing $\Omega$, we consider a corresponding
discretization (i.e., nodal representation) of $m$ and $u$, as well as the differential operator.
%Upon discretizing $m$ and $f$ on a given mesh and 
Reshaping these mesh functions into vectors we can write the resulting approximation
of the forward operator as
\begin{eqnarray}
\ff (\mm) =  P \uu =  P G(\mm) \qq,
\label{2.8}
\end{eqnarray}
where the right hand side vector $\qq$ is commonly referred to as
a source, $G$ is a discrete Green's function (i.e., the inverse of the discretized PDE operator), 
$\uu = G(\mm)\qq$ is the field (i.e., the PDE solution, here an interim quantity),
and $P$ is a projection matrix that projects the field to the locations where the data values $\dd$ 
are given. 

This setup is typical in the thriving research area of inverse problems; 
see, e.g.,~\cite{ehn1,vogelbook}. A specific example is provided in Section~\ref{sec:application}.

The inverse problem is to find {$\mm$} such that the predicted and observed
data agree to within noise $\beata$: ideally,
\begin{eqnarray}
\dd  = \ff(\mm) + \beata . \label{2.9}
\end{eqnarray}
To obtain such a model $\mm$ that satisfies~\eqref{2.9} we need to estimate
the {\em misfit function} $\phi(\mm)$, i.e., the normed difference between observed data $\dd$ and predicted data $\ff(\mm)$.
An iterative algorithm is then designed to sufficiently reduce this misfit function.
But, {\em which norm should we use to define the misfit function?}

It is customary to conveniently assume that the noise satisfies 
$\beata \sim \mathcal{N}(0,\sigma I)$, i.e., that the noise is normally distributed with
a scaled identity for the covariance matrix, where $\sigma$ is the
standard deviation.
Then the  maximum likelihood (ML) {data misfit} function is simply the squared $\ell_2$-norm\footnote{
For a more general symmetric positive definite covariance matrix $\Sigma$, 
such that $\beata \sim \mathcal{N}(0,\Sigma)$,
we get weighted least squares, or an ``energy norm'', with the weight matrix $\Sigma^{-1}$ for $\phi$.
But let's not go there in this article.}
\begin{equation}
		\phi(\mm) = \| \ff (\mm) - \dd \|_2^{2} .
		\label{2.10}
		\end{equation}
In this case,
the celebrated Morozov {\em discrepancy principle} yields the stopping criterion 
\begin{eqnarray}
  \phi(\mm) \leq \rho , \quad {\rm where~} \rho = \sigma^2 l ,
	\label{2.11}
	\end{eqnarray}
see, e.g.,~\cite{ehn1,morozov,kaso}.	
So, here is a class of problems where we do have a meaningful and directly usable
tolerance value!

Assuming that a known tolerance $\rho$ must be satisfied as in \eqref{2.11}
is often too rigid in practice, because realistic data do not quite satisfy the assumptions
that have led to~\eqref{2.11} and~\eqref{2.10}.
Well-known techniques such as L-curve and GCV (see, e.g.,~\cite{vogelbook}) are specifically designed to handle more general
and practical cases where~\eqref{2.11} cannot be used or justified.
Also, if~\eqref{2.11} is used then a typical algorithm would try to find $\mm$
such that $\phi (\mm)$ is (smaller but) not much smaller than $\rho$,
because having $\phi (\mm)$ too small would correspond
to fitting the noise -- an effect one wants to avoid. 
The latter argument and practice do not follow from~\eqref{2.11}.

Moreover, keeping the misfit function $\phi(\mm)$ in check does not necessarily
imply a quality reconstruction 
(i.e., an acceptable approximation $\mm$ for the ``true solution'' $\mm_{\rm true}$,
which can be an elusive notion in itself). However, $\phi (\mm)$, 
and not direct approximations of $\| \mm_{\rm true} - \mm \|$,
is what one typically has to work with.\footnote{
The situation here is different from that in Section~\ref{sec2.1}, where the choice of local
error criterion over a global one was made based on convenience and efficiency considerations.
Here, although controlling $\phi(\mm)$  is merely a necessary and not sufficient condition
for obtaining a quality reconstruction $\mm$, it is usually all we have to work with.}
 So any additional a priori information
is often incorporated through some regularization.
  
Still, despite all the cautious comments in the preceding two paragraphs,
we have in~\eqref{2.11} in a sense a more meaningful practical expression 
for stopping an iterative algorithm than hitherto.

%%%%%%%%%%%%%%%%%%%%%%%%%%%%%%%%%%%%%%%%%%%%%%%%%%%%%%%

\subsubsection{Regularization and constrained formulations}
\label{sec2.3.1}

Typically there is a need to regularize the inverse problem, and often this is
done by adding a regularization term to~\eqref{2.10}. Thus, one attempts to 
{\em approximately} solve the Tikhonov-type problem
\begin{eqnarray}
\min_{\mm}~ \phi (\mm) + \lambda R(\mm), \label{2.12}
\end{eqnarray}
where $R(\mm) \geq 0 $ is a prior (we are thinking of some norm or semi-norm of $\mm$),
and $\lambda \geq 0$ is a regularization parameter.
Two other forms of the same problem~\eqref{2.12} immediately arise upon
interpreting $\lambda$ or $\lambda^{-1}$ as a Lagrange multiplier. These are the constrained
optimization formulations
\begin{eqnarray}
\min_\mm~ \phi (\mm), \quad {\rm s.t.~} R(\mm) \leq \tau; \ \ {\rm and} \label{2.13}
\end{eqnarray}
%and
\begin{eqnarray}
\min_\mm~  R(\mm) , \quad {\rm s.t.~}  \phi (\mm) \leq \rho . \label{2.14}
\end{eqnarray}
The non-negative parameters $\rho, \ \tau$ and $\lambda$ are related to one another
in a nontrivial manner that makes these three formulations indeed equivalent; see, e.g.,~\cite{befr08}.

Each of these formulations has its fan club. The inevitable discussion regarding which is
best becomes more heated when $R$ involves the $\ell_1$-norm, corresponding to a prior that
favours some form of sparsity; see~\cite{doasha} and references therein.
This is so not only because sparsity is popular and extremely useful, 
but also because use of the $\ell_1$-norm introduces
lower smoothness in $R(\mm)$.
In such circumstances, the formulation~\eqref{2.13} can be more directly amenable
to efficient solution techniques (see the equivocal~\cite{befr08}), while~\eqref{2.14} has the advantage that
the tolerance $\rho$ is known (cf.~\eqref{2.11}) with far less a priori uncertainty than $\tau$ or $\lambda$. 
Thus, the popular question of choosing the most practically advantageous formulation from among
\eqref{2.12}--\eqref{2.14} (see, e.g.,~\cite{hatifr})
is tied to the topic of the present article.

%%%%%%%%%%%%%%%%%%%%%%%%%%%%%%%%%%%%%%%%%%%%%%%%%%%%%%%

\subsection{Least squares data fitting with many data sets}
\label{sec2.4}

%Let us next 
In our fourth case study we generalize the setting of Section~\ref{sec2.3}, allowing
not only one but several (indeed, many) data sets $\dd_i, \ i = 1, \ldots , s$,
where each data set has length $l$, so $\dd_i \in {\R}^{l}$.
Unless $s$ is small, working with such a problem can be prohibitively expensive,
so we next use randomized approximations to the misfit function.
This in turn provides a new meaning to the uncertainty in the given tolerance $\rho$, which we further explore
in later sections.

We can conveniently define an $l \times s$ data matrix $D$ whose
$i$th column is $\dd_i$.
Correspondingly, there are forward operators
\begin{subequations}
\begin{eqnarray}
\ff_i (\mm) &\equiv&  \ff (\mm,\qq_i) =  P G(\mm) \qq_i, \quad i = 1, \ldots, s,
\label{3.12a}
\end{eqnarray}
where $\qq_i$ are different sources, and we
arrange the forward operator in form of an $l \times s$ matrix function $F (\mm)$
which has $\ff_i$ as its $i$th column.
The inverse problem is to find {$\mm$} such that the predicted and observed 
data agree to within noise,  
\begin{eqnarray}
\dd_i  = \ff_i(\mm) + \beata_i , \quad i = 1, 2, \ldots, s. \label{3.12b}
\end{eqnarray}
Assuming next that
$\beata_i \sim \mathcal{N}(0,\sigma I)$,
the ML data misfit function is
\begin{equation}
		\phi(\mm) = \sum_{i=1}^s\| \ff_i (\mm) - \dd_i \|_2^{2} = \| F(\mm) - D \|_F^2 ,
		\label{3.12c}
		\end{equation}
where the subscript $F$ denotes the Frobenius norm. In this case,
the discrepancy principle yields the stopping criterion 
\begin{eqnarray}
  \phi(\mm) \leq \rho , \quad {\rm where~} \rho = \sigma^2 ls .
	\label{3.12d}
	\end{eqnarray}
	So, as in Section~\ref{sec2.3}, provided that all assumptions hold
	we have a decent idea of a stopping tolerance value
	for an iterative method that generates iterates $\mm_k, k = 1, 2, \ldots$, in an attempt
	to decrease $\phi$ until~\eqref{3.12d} holds with approximate equality.
	%(We note in passing that some extensions of the probability distribution assumption on the noise,
	%which lead to {\em weighted} least squares generalizing~\eqref{3.12c}, are possible; see~\cite{roszas}.)
	\label{3.12}
	\end{subequations}
	
However, when $s$ is very large, since $s$ evaluations of $G(\mm) \qq_i$ are required just to form $F(\mm)$ for a given
$\mm$,  we obtain an instance where the matrix $B = F(\mm) - D$ can be
prohibitively expensive to calculate explicitly (Section~\ref{sec:application} provides a concrete example).
Hence the matrix $A = B^TB$ is implicit symmetric positive semi-definite (SPSD) (cf. Section~\ref{sec2.2}).
Furthermore, we have
 \begin{equation}
\phi(\mm) = \| B \|_{F}^{2} = tr (A) = \Ex(\| B \ww \|_{2}^{2}),
\label{frob_trace}
\end{equation}
where $tr(A)$ denotes the trace of the matrix $A$, $\Ex$ stands for expectation,
and $\ww$ stands for a random vector drawn from any distribution satisfying $\Ex(\ww \ww^{T}) = \mathbb{I}$. 
%($\mathbb{I}$ is the identity matrix) 
%$tr(A)$ denotes the trace of the matrix $A$, and $\Ex$ denotes the expectation.

Inexpensive approximate alternatives to working with all $s$ data sets throughout
an algorithm for finding a suitable $\mm$ may be obtained by   
approximating the misfit function~$\phi(\mm)$ in~\eqref{3.12c}.
Specifically, the rightmost expression in~\eqref{frob_trace}
%By~\eqref{frob_trace} this is equivalent to approximating the corresponding matrix trace.
%and the Monte Carlo methods discussed in Section~\ref{sec:trace} apply. For the latter, 
suggests a Monte Carlo sampling method, obtaining
% may be devised and analyzed.
%
%By~\eqref{3.3} we obtain
\begin{equation}
\widehat{\phi} (\mm,n) \defeq \frac{1}{n} \sum_{j=1}^{n} \|B(\mm) \ww_{j}\|_{2}^{2} \approx \phi(\mm).
\label{approx_phi}
\end{equation}
Note that $\widehat{\phi}(\mm,n)$ is an {\em unbiased estimator} of $\phi(\mm)$, as we have 
%\begin{equation*}
$\phi(\mm) = \Ex (\widehat{\phi}(\mm,n))$.
For the forward operators \eqref{3.12a}, if $n \ll s$ (in words, the sample size $n$ is much smaller
than the total number of data sets $s$) then 
this procedure yields a very efficient algorithm for approximating the misfit~\eqref{3.12c}. 
This is so because 
%(letting $w_{i}$ denote the $i^{th}$ component of a realization of random vector $\ww$) 
\begin{equation*}
\sum_{i=1}^{s} \ff(\mm,\qq_{i}) w_{i}  = \ff(\mm,\sum_{i=1}^{s} w_{i} \qq_{i} ) ,
\end{equation*}
which can be computed with a single evaluation of $\ff$ per realization of the random vector 
$\ww = (w_1, \ldots , w_s)^T$, so in total $n$ (rather than $s$) such evaluations 
are required~\cite{HaberChungHermann2010,rodoas1}.
An example utilizing~\eqref{approx_phi} is considered in Section~\ref{sec:application} below.

But now it is natural to ask, {\em how large should $n$ be?}
Generally speaking, the larger $n$ is, the closer $\widehat{\phi}(\mm,n)$
is to $\phi(\mm)$. So we are led to ask how important it is to pedantically satisfy~\eqref{3.12d}.
This in turn brings up the question that is the common refrain of this section,
namely, the extent to which the stopping criterion is really known and as such must be obeyed.  
Here, a higher uncertainty in $\rho$ and~\eqref{3.12} allows a more efficient solution algorithm
because a smaller $n$ will suffice.  
%The answer to this 
%we also have a new uncertainty associated with the tolerance $\rho$,
%because it appears in~\eqref{3.12d} as a bound on $\phi(\mm)$, not $\widehat{\phi} (\mm,n)$.
In the context of the present randomized algorithm it is therefore natural to consider
satisfying the error criterion~\eqref{1} only within some probability range.
To concentrate on this aspect we isolate it by assuming that $\rho$
in~\eqref{3.12d} is given and the error model that enables this is approximately valid.

%%%%%%%%%%%%%%%%%%%%%%%%%%%%%%%%%%%%%%%%%%%%%%%%%%%%%%%%%%%%%%%%%%%%%%%%%

%\section{Satisfying a stopping criterion, probably} 
\section{Probabilistic relaxation of a stopping criterion} 
\label{sec:prob}

The previous section details four different case studies which highlight the fact
of life that in applications an error tolerance for stopping an algorithm is rarely known with
absolute certainty. Thus, we can say that such a tolerance is only ``probably'' known. 
%where an error tolerance (or its meaning) can be relaxed
%in various ways in order to achieve some pragmatic computational or algorithmic advantage.
%In Section~\ref{sec2.1} this is expressed by moving from global to local error estimate,
%hopefully cautiously. In Section~\ref{sec2.2} we point out that the value of a tolerance
%may not have much to do with application, yet some caution is required there as well.
%In Section~\ref{sec2.3} we see for the first time some realistic tolerance value,
%obtained albeit under some potentially restrictive assumptions.
%In section~\ref{sec2.4} we see for the first time some probabilistic issues arising
%in connection with the meaning of a given tolerance value.
%The word ``probably'' can be considered in different ways. One is an everyday meaning,
%where we say ``probably'', meaning a mixture of ``I'm not sure'' and ``I don't care''.
Yet in some situations, it is also possible to assign it a more precise meaning
in terms of statistical probability. This holds true for the case study in Section~\ref{sec2.4},
with which we proceed below.
%Hopefully, at this point we have made our point that considering an error tolerance as
%determinstically given with 100\% assurance is usually not reasitic, so we are permitted
%to embark on the following probabilistic foray.   
Thus, in the present section we consider a way to relax~\eqref{1},
which is more systematic and also allows for further theoretical developments.
Specifically, we consider
%relax the concept to 
satisfying a tolerance in a probabilistic sense. 

Suppose we seek an unbiased estimator $\hat g(x)$ to a real 
%(complex) 
valued function $g(x)$
(so $\Ex (\hat g(x)) = g(x)$, where $\mathbb{E}$ denotes expectation).
 Then, given a pair of values {$(\veps,\delta)$}, both small and positive, we require
		\begin{equation}
			Pr \big(| \hat g(x) - g(x) | \leq \veps~|g(x)| \big) \geq 1-\delta .
			\label{3.1}
		\end{equation}
The parameters $\veps$ and $\delta$ relate to the relative accuracy and the probabilistic 
guarantee of such an estimation, respectively~\cite{achlioptas}.
%This idea is well-known in some but not all communities.
In Section~\ref{sec:trace} we consider applying this notion to the problem of trace estimation.
Then in Section~\ref{sec:nls} we apply the results of Section~\ref{sec:trace}
to the case study of Section~\ref{sec2.4}.
The uncertainty in the stopping tolerance $\rho$ is quantified in terms of parameters
$\veps$ and $\delta$ towards the end of this section.
% and methods where such a notion becomes handy.

%%%%%%%%%%%%%%%%%%%%%%%%%%%%%%%%%%%%%%%%%%%%%%%%%%%%%%%%%%%%%%%%%%%%%%%%

\subsection{Estimating the trace of an implicit matrix}
\label{sec:trace}

As in Section~\ref{sec2.2} we consider an $s \times s$ matrix $A$ that is given only
implicitly, through matrix-vector products. We assume that $A = B^TB$,
%symmetric positive semi-definite (SPSD), 
where $B$ is some real rectangular matrix
with $s$ columns, and wish to approximate the
trace, $tr(A) = \sum_{i=1}^s a_{i,i}$, without having the luxury of knowing
the diagonal elements $a_{i,i}$. 
This task has several applications; see~\cite{avto,roas1}.
We consider it here mainly as a preparatory step, given the appearance
of the trace in~\eqref{frob_trace}.
%and in Sections~\ref{sec:nls} and~\ref{sec:application}
%we consider one such; see~\cite{avto,roas1}.

Note that if $\xx$ is the $i$th column of the $s \times s$
identity matrix then $(\xx ,A\xx) = a_{i,i}$. Hence, $tr(A)$ can be calculated
in $s$ matrix-vector products. However, $s$ is very large so we want a cheaper
approximation.
Towards that goal, observe that if $\ww$ stands for a random vector drawn
from a probability distribution ${\cal D}$ satisfying
$\mathbb{E} [ \ww \ww^T ] = I$,  then
\begin{eqnarray}
tr(A) = \Ex [(\ww, A \ww) ] .
\label{3.2}
\end{eqnarray}
So, we consider the Monte Carlo approximation $tr(A) \approx tr_{\cal D} (A)$, 
%with the approximation 
defined by
\begin{eqnarray} 
tr_{\cal D} (A) \defeq \frac{1}{n} \sum_{i=1}^{n} (\ww_{i}, A \ww_{i}) = 
\frac{1}{n} \sum_{i=1}^{n} \| B\ww_i \|_2^2 \ ,
\label{3.3}
\end{eqnarray}
where $\ww_i$ are drawn from the distribution ${\cal D}$, and hopefully $n \ll s$.

The next question is, {\em how small can we take $n$ to be?}
In other words, how can we quantify the uncertainty in such an approximation?
Here is where the pair {$(\veps,\delta)$} appearing in~\eqref{3.1}
comes in handy for deriving probabilistic necessary and sufficient conditions on the size of $n$.
%in an attempt to quantify the uncertainty involved.
Below we state two 
%some of the 
results that were proved in~\cite{roas1,roszas}.

Let us restate~\eqref{3.1} in the present context as
\begin{equation}
			Pr \big(| tr_{\cal D} (A) - tr(A) | \leq \veps~|tr(A)| \big) \geq 1-\delta .
			\label{3.4}
		\end{equation}
Further, given the pair $(\veps,\delta)$ (both values positive and small), define the constant
\begin{eqnarray}
c = c( \veps,\delta) = \veps^{-2} \ln (2/\delta).
\label{3.5}
\end{eqnarray}
This constant appears in all estimates of the type stated in Theorem~\ref{thm1}
below. Note its strong
dependence on $\veps$ and much weaker dependence on	$\delta$, in the sense of how rapidly
$c$ grows as these values shrink.	

\begin{theorem}
$  $
\begin{itemize}
\item
Let ${\cal D}$ be the Gaussian (i.e., standard normal) probability distribution. 
Then the probabilistic bound~\eqref{3.4} holds if
\[ n \geq 8 c (\veps,\delta) .\]
\item
Let ${\cal D}$ be the Rademacher probability distribution~\cite{hutchinson},
namely, for each component of $\ww = (w_1, \ldots , w_s)^T$, $Pr(w_j = 1) = Pr(w_j = -1) = 1/2$.
(The components of $\ww$ are i.i.d.)
 Then the probabilistic bound~\eqref{3.4} holds if
\[ n \geq 6 c (\veps,\delta) .\]
\end{itemize}
\label{thm1}
\end{theorem}

The property that stands out in these bounds is that they are independent of the matrix size $s$
(reminiscent in this sense to a multigrid method for the Poisson example in Section~\ref{sec2.2})
and of any property of the matrix $A$ other than it being SPSD.
The latter also has a downside, of course, in suggesting that such bounds may not be very tight.
%Indeed, if $A$ happens to be a diagonal matrix (although we pretend to not know that) then using
%the Rademacher distribution 
%we obtain the trace precisely with $n = 1$, whereas with the Gaussian distribution $n$ may be quite large.
%Such a difference is not reflected in the bounds of Theorem~\ref{thm1}. In~\cite{roas1}
%there are additional theorems that explore satisfying bounds of the form~\eqref{3.4},
%taking into account matrix properties that may be known to hold even if they cannot be easily
%verified.

Theorem~\ref{thm1} provides {\em sufficient} bounds on $n$, and not very sharp ones at that;
however, these bounds do have a charming simplicity. 
Next we give better sufficient bounds, as well as {\em necessary} bounds,
for the case where ${\cal D}$ is the Gaussian distribution~\cite{roszas}.
Simplicity, however, shall have to be sacrificed.
For this purpose, we write~\eqref{3.4}, with a minor abuse of notation, as
\begin{subequations}
\begin{eqnarray}
\Pr\Big( tr_{n}(A) &\geq& (1-\veps) tr(A) \Big) \geq 1-\delta \label{prob_ineq_lower}, \\
\Pr\Big( tr_{n}(A) &\leq& (1+\veps) tr(A) \Big) \geq 1-\delta \label{prob_ineq_upper}.
\end{eqnarray}
%Recall further that the gamma distribution parametrized by shape $\alpha$ and rate $\beta$ %parameters 
%has the probability density function
%\begin{eqnarray*} 
%f(x) = \begin{cases}
%\frac{\beta^{\alpha}}{\Gamma(\alpha)} x^{\alpha - 1 } e^{-\beta x}  &\text{$x \ge 0$} , \cr
%0 &\text{$x < 0$} .\end{cases} \label{gamma_pdf}
%\end{eqnarray*}
Also, denote by
%$n_{0}$ and $n$ denote positive integers, and 
$Q_n \sim \chi^{2}_{n}$ 
a chi-squared random variable of degree $n$,
%(which is a special case of the gamma distribution), 
and set 
$ Q(n) = \frac {Q_n}n$.
% We have~\cite{roszas}
\label{prob_ineq_lower_upper}
\end{subequations}

%We have
\begin{theorem}[Necessary and sufficient condition for~\eqref{prob_ineq_lower_upper}]
Given an SPSD matrix $A$ of rank $r$ and parameters $(\veps,\delta)$ as above, the following hold:
\begin{enumerate}[(i)]

\item {\em Sufficient} condition for~\eqref{prob_ineq_lower}: there exists some integer $n_{0} \geq 1$ such that
	\begin{equation}
			\Pr\left( Q(n_{0}) < (1-\veps) \right) \leq \delta .
		\label{max_bnd_N_lower}
	\end{equation}
Furthermore,~\eqref{prob_ineq_lower} holds for all $n \geq n_{0}$.

		\item {\em Sufficient} condition for~\eqref{prob_ineq_upper}: if the inequality 
	\begin{equation}
			\Pr\left( Q(n_{0}) \leq (1+\veps) \right) \geq 1-\delta
		\label{max_bnd_N_upper}
	\end{equation}
	is satisfied for some $n_{0} > \veps^{-1}$, then~\eqref{prob_ineq_upper} holds with $n = n_{0}$. 
	Furthermore, there is always an $n_{0} > \veps^{-2}$ such that \eqref{max_bnd_N_upper}
	is satisfied and, for such $n_{0}$, it follows that~\eqref{prob_ineq_upper} holds for all $n \geq n_{0}$.
	
	\item {\em Necessary} condition for~\eqref{prob_ineq_lower}: if~\eqref{prob_ineq_lower} holds for some $n_{0} \geq 1$, 
	then for all $n \geq n_{0}$
	\begin{equation}
			\Pr\left( Q(n r) < (1-\veps)\right) \leq \delta.
		\label{min_bnd_N_lower}
	\end{equation}

	%\item \textbf{Existence of such $\bm {N_{0}}$:} for any pair $(\veps,\delta)$, there exists $N_{0} > \frac{1}{\veps^{2}}$ such that~\eqref{max_bnd_N_upper} holds,
		%\item \textbf{Tightness:} if $A$ has rank one, there is a smallest $N_{0} > \frac{1}{\veps^{2}}$ satisfying~\eqref{max_bnd_N_upper} such that for any $N \geq N_{0}$,~\eqref{prob_ineq_upper} holds, and for all $\frac{1}{\veps^{2}} < N < N_{0}$,~\eqref{prob_ineq_upper} does not hold. If $\delta$ is small enough so that~\eqref{prob_ineq_upper} does not hold for any $ N \leq \frac{1}{\veps^{2}}$, then $N_{0}$ is necessary and sufficient for~\eqref{prob_ineq_upper},
	\item {\em Necessary} condition for~\eqref{prob_ineq_upper}: if~\eqref{prob_ineq_upper} 
	holds for some $n_{0} > \veps^{-1}$, then
	\begin{equation}
			\Pr\left(Q(n r) \leq (1+\veps)\right) \geq 1- \delta,
		\label{min_bnd_N_upper}
	\end{equation}
	with $n = n_{0}$. Furthermore, if $n_{0} > \veps^{-2} r^{-2}$, 
	then~\eqref{min_bnd_N_upper} holds for all $n \geq n_{0}$.
	
	\end{enumerate}
	\label{main_trace_theorem}
\end{theorem}

%%%%%%%%%%%%%%%%%%%%%%%%%%%%%%%%%%%%%%%%%%%%%%%%%%%%%%%%%%%%%%%%%%%%%%%%%%%%%%%

The necessary conditions in Theorem~\ref{main_trace_theorem} indicate that the lower bound on the {\em smallest} 
``true'' $n$ that satisfies~\eqref{3.4} grows as the rank of $A$ decreases (regardless of $A$'s size $s, \ s \geq r$).
For $r = 1$ the necessary and sufficient bounds coincide, which indicates a form of tightness.

{\em Are these necessary bounds ``a bug or a feature''?}
We argue that they can be both. Examples can be easily found where Hutchinson's method~\cite{hutchinson}
(employing Rademacher's distribution) performs better
than Gaussian, requiring a smaller $n$ for a small rank matrix $A$.
On the other hand, if other factors come in (as in the application of Section~\ref{sec:application}
below) which make the practical use of the Gaussian and Hutchinson methods comparable, then the availability
of both necessary and sufficient conditions for the Gaussian distribution allow a better
chance to quantify the error, probabilistically, in cases where the lower and upper bounds are close.
This is taken up next.
%below, following~\cite{roszas}. 
%see~\cite{}

%%%%%%%%%%%%%%%%%%%%%%%%%%%%%%%%%%%%%%%%%%%%%%%%%%%%%%%%%%%%%%%%%%%%%%%%

\subsection{Least squares data fitting with many data sets revisited}
\label{sec:nls}

We now return to the problem considered in Section~\ref{sec2.4}, and apply the results
of the previous section to~\eqref{approx_phi} and~\eqref{3.12d}.
Thus, according to~\eqref{prob_ineq_lower_upper},
in the check for termination of our iterative algorithm at the next iterate $\mm_{k+1}$, 
we consider replacing the condition
\begin{subequations}
\begin{eqnarray}
\phi (\mm_{k+1}) \leq \rho \label{3.13a}
\end{eqnarray}
by either
\begin{eqnarray}
\widehat{\phi}(\mm_{k+1},n_{t}) &\leq& (1-\veps) \rho , \quad {\rm or} \label{stop_crit_hard} \\
\widehat{\phi}(\mm_{k+1},n_{t}) &\leq& (1+\veps) \rho , \label{stop_crit_soft}
\end{eqnarray}
for a suitable $n = n_t$ that is governed  by Theorem~\ref{main_trace_theorem} with a prescribed pair $(\veps,\delta)$. 
If~\eqref{stop_crit_hard} holds, then it follows with a probability of 
at least  $(1-\delta)$ that~\eqref{3.13a} holds. 
On the other hand, if~\eqref{stop_crit_soft} does {\em not} hold, 
then we can conclude with a probability of at least  $(1-\delta)$ that~\eqref{3.13a} is {\em not} satisfied. 
In other words, unlike~\eqref{stop_crit_hard}, a successful~\eqref{stop_crit_soft} is only necessary and not 
sufficient for concluding that~\eqref{3.13a} 
holds with the prescribed probability $1-\delta$.
%for a suitable $n = n_t$ that is governed  by Theorem~\ref{main_trace_theorem}.
%Note that a 
%successful~\eqref{stop_crit_soft} is only necessary and not sufficient for concluding that~\eqref{3.13a} 
%holds with the prescribed probability $1-\delta$.
\label{stop_crits}
\end{subequations}

{\em What are the connections among these three parameters, $\rho, \ \delta$ and $\veps$?!}
The parameter $\rho$ is the deterministic but not necessarily too trustworthy error tolerance
appearing in~\eqref{3.13a},
much like the tolerance in Section~\ref{sec2.1}.
Next, we can reflect the uncertainty in the value of $\rho$ by choosing an appropriately large $\delta ~(\leq 1)$. 
Smaller values of $\delta$ reflect a higher certainty in $\rho$ and a more rigid stopping criterion
(translating into using a larger $n_t$).
For instance, success of~\eqref{stop_crit_hard} is equivalent to making a statement 
on the probability that a positive ``test'' result will be a ``true'' positive. 
This is formally given by the conditional probability statement 
%For instance, using~\eqref{stop_crit_hard}, we are making 
%a statement on conditional probability:
\[ Pr \left( \phi (\mm_{k+1} ) \leq \rho ~|~ \widehat{\phi}(\mm_{k+1} , n_t ) \leq (1 - \veps)\rho\right)Ê \geq 1 - \delta. \]
Note that, once the condition in this statement is given, 
the rest only involves $\rho$ and $\delta$. So the tolerance $\rho$ is augmented by the probability parameter $\delta$.
The third parameter $\veps$ governs the false positives/negatives 
(i.e., the probability that the test will yield a positive/negative result, if in fact~\eqref{3.13a} is false/true), 
where a false positive is given by
\begin{eqnarray*}Ê
Pr \left( \widehat{\phi}(\mm_{k+1} , n_t ) \leq (1 - \veps)\rho ~|~ \phi (\mm_{k+1} ) > \rhoÊ \right) ,
\end{eqnarray*}
while a false negative is
\begin{eqnarray*}
Pr \left( \widehat{\phi}(\mm_{k+1} , n_t ) > (1 - \veps)\rho ~|~ \phi(\mm_{k+1} ) \leq \rhoÊ\right) . %\blacklozenge
\end{eqnarray*} 

%{\bf Remark:}	
We note in passing that such efficient algorithms as described above
can also be obtained with some extensions of the probability distribution assumption on the noise, 
which lead to {\em weighted} least squares generalizing~\eqref{3.12c}.
%; see~\cite{roszas}.

%some extensions of the probability distribution assumption on the noise,
	%which lead to {\em weighted} least squares generalizing~\eqref{3.12c}, are possible; see~\cite{roszas}.
	
%The upshot of this section is that 
In summary, we have demonstrated in this section an approach for quantifying
the	uncertainty in the stopping criterion
%an error tolerance $\rho$ 
by deriving such a procedure for the case study described
in Section~\ref{sec2.4}.

%%%%%%%%%%%%%%%%%%%%%%%%%%%%%%%%%%%%%%%%%%%%%%%%%%%%%%%%%%%%%%%%%%%%%%%%%

\section{An inverse problem with many data sets} 
\label{sec:application}

The purpose of this section is to demonstrate the ideas developed in Section~\ref{sec:prob}
for the methods of Section~\ref{sec2.4} on a concrete example.

Inverse problems of the sort described in~\eqref{3.12} arise frequently
in practice, and applications include electromagnetic data inversion in mining exploration 
(e.g.,~\cite{na,dmr,haasol,olhash}), 
seismic data inversion in oil exploration (e.g.,~\cite{fichtner,hel,rnkkda}), 
diffuse optical tomography (DOT) (e.g.,~\cite{arridge,boas}), 
quantitative photo-acoustic tomography (QPAT) (e.g.,~\cite{gaooscher,yuan}), 
direct current (DC) resistivity (e.g.,~\cite{smvoz,pihakn,haheas,HaberChungHermann2010,doas3}),  
and electrical impedance tomography (EIT) (e.g.,~\cite{bbp,cin,doasha}).
Exploiting many data sets currently appears to be particularly popular in exploration geophysics,
and our example can be viewed as mimicking a DC resistivity setup.

Our entire setup is the same as in~\cite{roszas} (built in turn on~\cite{doas3,rodoas1}), where many details that are
omitted here can be found.
% see also~\cite{doas3,rodoas1}.
%because this is not what we 
This allows us to concentrate next on those issues that are most relevant in the present article.
The PDE has the form
\begin{subequations}
\begin{eqnarray}
\div (\mu(\xx) \grad u) = q(\xx), \quad \xx \in  \Omega ,
\label{42.1a}
\end{eqnarray}
where $\Omega \subset \R^d$, $d = 2$ or $3$, and
$\mu(\xx) \geq \mu_0 > 0$ is a conductivity function which may be rough 
(e.g., only piecewise continuous).
An appropriate transfer function $\psi$ is selected so that, point-wise, $\mu(\xx) = \psi(m(\xx))$. 
For example, $\psi$ can be chosen so as to
ensure that the conductivity stays positive and bounded away from $0$, as well as to incorporate bounds,
which are often known in practice.
The matrix $G$ in~\eqref{3.12a} is a discrete Green's function for \eqref{42.1a}
subject to the homogeneous Neumann boundary conditions
\begin{eqnarray}
\frac{\partial u}{\partial n} = 0, \quad \xx \in \partial\Omega. \label{42.1b}
\end{eqnarray}
In other words, evaluating $\ff_i (\mm)$ for given $\mm$ and $i$ requires the approximate
solution of the PDE problem~\eqref{42.1}.
Given the nonlinearity in $\mm$ (which requires an iterative method for the nonlinear least squares
problem of minimizing $\phi (\mm)$) and the number of data sets $s$, 
the PDE count can easily climb without employing the Monte Carlo
approximations described in Sections~\ref{sec2.4} and~\ref{sec:nls}.
%For $\Omega$ we will consider a unit square or a unit cube.
%In what follows, we will consider sources $\qq$ that are sums of $\delta$-functions.
\label{42.1}
\end{subequations}

The variant of our algorithm 
%presented in~\cite{roszas} and 
used here
employs unbiased estimators of the form~\eqref{approx_phi}
on four occasions during an iteration $k$: 
\begin{enumerate}
\item
Use sample size $n_k$ for an approximate stabilized Gauss-Newton (ASGN) iteration. This requires
also a corresponding approximate gradient, and $n_k$ is set heuristically
from one iteration to the next as described next.
\item
Use sample size $n_c$ for a cross validation step, where we check whether
\begin{equation}
(1 - \veps) \widehat{\phi}(\mm_{k+1},n_{c}) \leq (1+\veps)  \widehat{\phi}(\mm_{k},n_{c})
\label{cross_valid_soft}
\end{equation}
holds. If it does not then we judge that the sample size $n_k$ is too small for a meaningful ASGN iteration,
increase it by a set factor (e.g., $2$) and reiterate; otherwise, continue. 
We set $n_c$ using a given pair $\veps_c, \delta_c$.
\item 
Use sample size $n_u$ for an uncertainty check, asking if
\begin{equation}
\widehat{\phi}(\mm_{k+1},n_{u}) \leq (1-\veps) \rho .
\label{uncert_check_hard}
\end{equation}
We set $n_u$ using a given pair $\veps_u, \delta_u$.
If~\eqref{uncert_check_hard} holds then we next check for possible termination.
\item
Use sample size $n_t$ for a stopping criterion check, asking if~\eqref{stop_crit_soft}
holds. If yes then terminate the algorithm.
We set $n_t$ using a given pair $\veps_t, \delta_t$.
\end{enumerate}  
The last three steps allow for uncertainty quantification according to Theorem~\ref{main_trace_theorem}.

\begin{figure}[htb]
\centering
\subfigure[]{\includegraphics[scale=0.25]{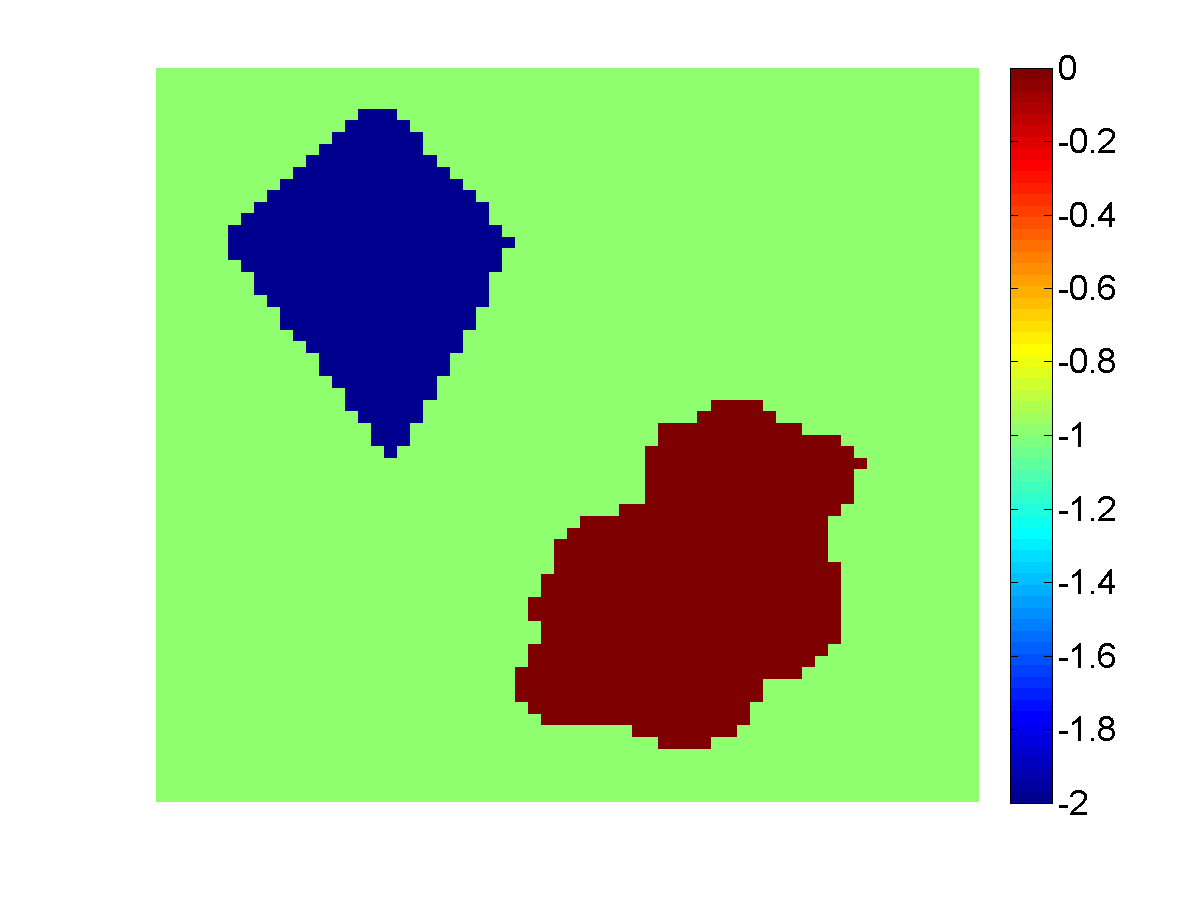}}
\subfigure[]{\includegraphics[scale=0.25]{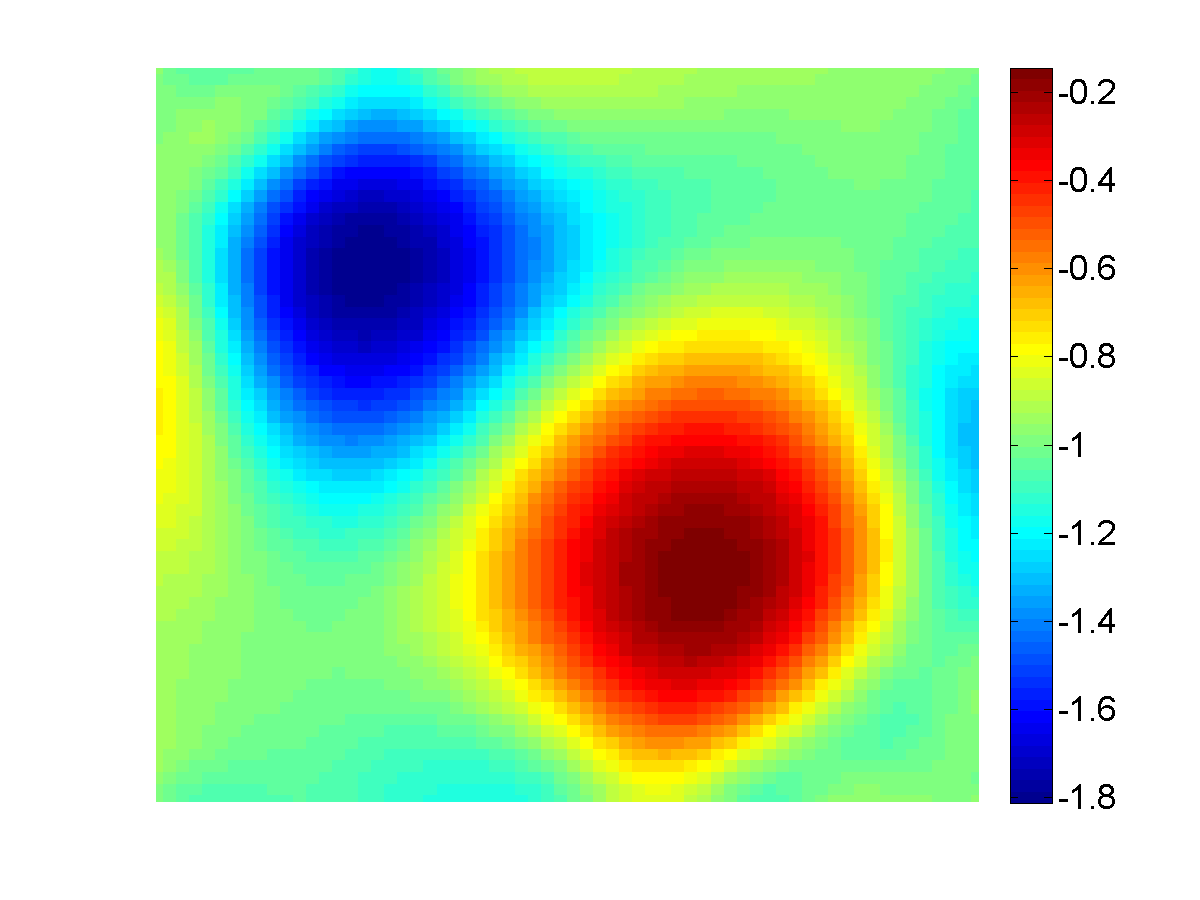}}
\subfigure[]{\includegraphics[scale=0.25]{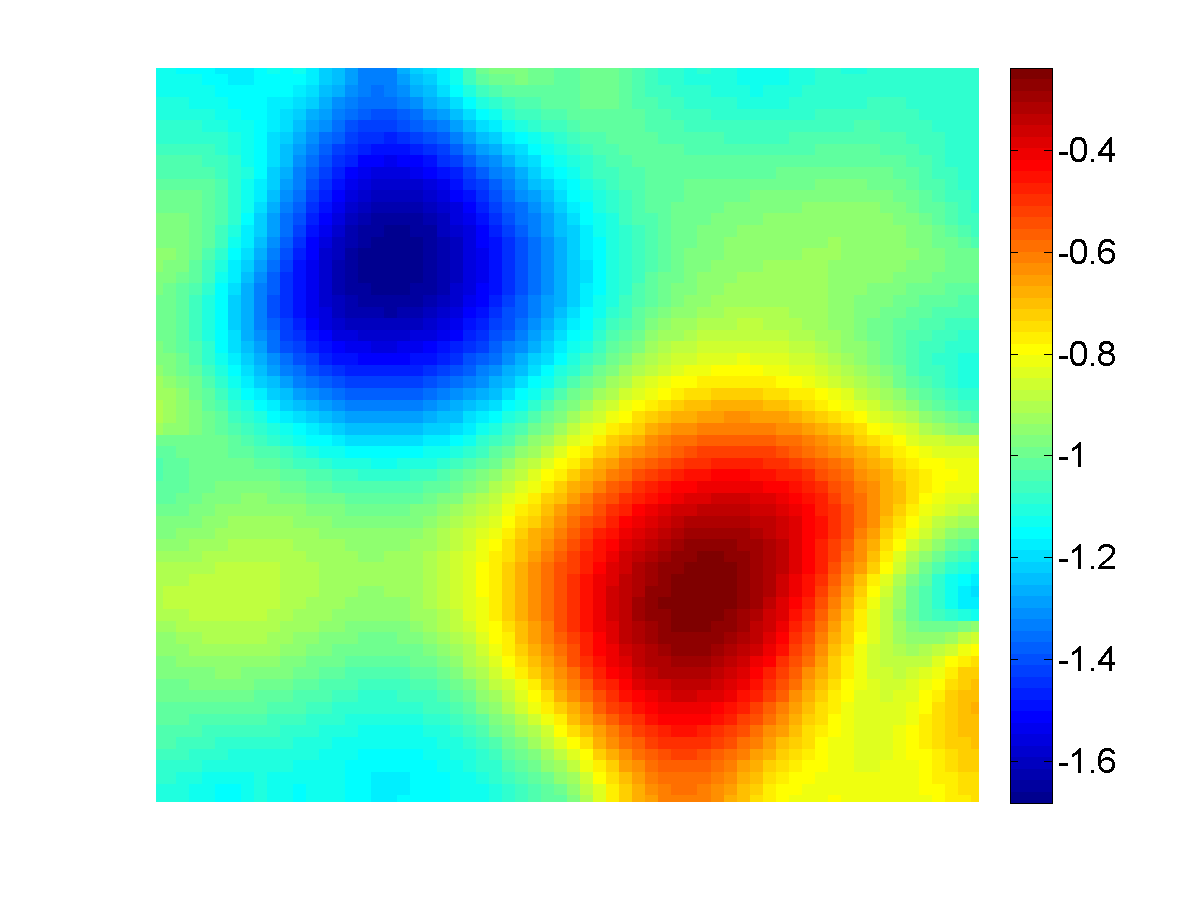}}
\subfigure[]{\includegraphics[scale=0.25]{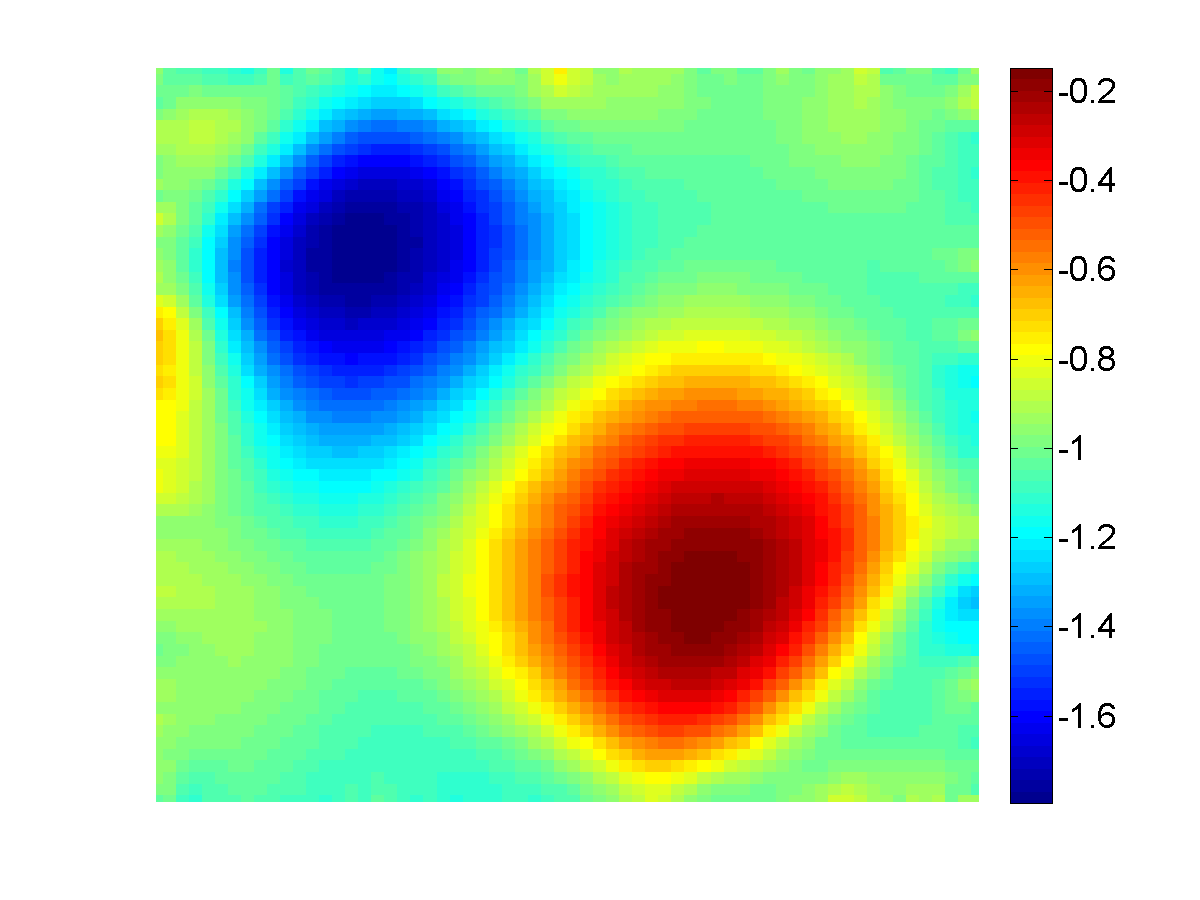}}
\caption{Plots of log-conductivity: (a) True model; (b) Vanilla recovery with $s = 3,969$; 
(c) Vanilla recovery with $s=49$; (d) Monte Carlo recovery with $s = 3,969$. 
The vanilla recovery using only $49$ measurement sets is clearly inferior, 
showing that a large number of measurement sets can be crucial for better reconstructions.
The recovery using our algorithm, however, is comparable in quality to Vanilla with
the same $s$. The quantifier values used in our algorithm were: $(\veps_c, \delta_c) = (0.05, 0.3)$,
$(\veps_u, \delta_u ) = (0.1, 0.3)$ and $(\veps_t, \delta_t) = (0.1, 0.1)$.
\label{fig02} }
\end{figure}
Figure~\ref{fig02} depicts an example of a piecewise constant ``true model'' $\mm_{\rm true}$ consisting of two
homogeneous but different bodies in a homogeneous background (a). This is used to synthesize noisy 
partial-boundary data sets that are then used in turn to approximately recover the model.
The large dynamical range of the conductivities, together with the fact that the data
is available only on less than half of the boundary, contribute to the difficulty in obtaining
good quality reconstructions.
The term ``Vanilla'' refers to using all $s$ available data sets for each task during the algorithm.
This costs 527,877 PDE solves\footnote{Fortunately, the discrete Green's function $G$ does not depend
on $i$ in~\eqref{3.12a}. Hence, if the problem is small enough that
a direct method can be used to construct $G$, i.e., perform one LU decomposition at each iteration $k$, 
then the task of solving half a million PDEs
just for comparison sake becomes less daunting.} for $s=3,969$ (b) and 5,733 PDE solves for $s=49$ (c). 
However, the quality
of reconstruction using the smaller number of data sets is clearly inferior.
On the other hand, using our algorithm 
%(see~\cite{roszas} for the full details)
yields a recovery (d) that is comparable to Vanilla but at the cost of only 
5,142 PDE solves. The latter cost is about 1\% that of Vanilla and is comparable in order of
magnitude to that of evaluating $\phi (\mm)$ once!

\subsection{TV and stochastic methods}
\label{sec4.1}

This section, like Section~\ref{sec2.2.1}, is not directly related to the main theme of
this article, but it arises from the present discussion and has significant merit
on its own. The two comments below are not strongly related to each other. 
\begin{itemize}
\item
The specific example considered above is used also in~\cite{roszas}, except that 
%the plots there are different because their 
the objective function there includes a total variation (TV) regularization.
This represents usage of additional a priori information (namely, that the true model is
discontinuous with otherwise flat regions), whereas here
an implicit $\ell_2$-based regularization has been employed without such knowledge regarding the true solution.
The results in Figures~4.3(b) and~4.4(vi) there correspond to our Figures~\ref{fig02}(b) and~\ref{fig02}(d),
respectively, and as expected, they look sharper in~\cite{roszas}.
On the other hand, a comparative glance at Figure~4.3(c) there vs the present Figure~\ref{fig02}(c) reveals
that the $\ell_1$-based technique can be inferior to the $\ell_2$-based one, even
for recovering a piecewise constant solution!
Essentially, even for this special solution form TV shines only with sufficiently good data,
and here ``sufficiently good'' translates to ``many data dets''.
This intuitively obvious observation does not appear to be as well-known today
as it used to be~\cite{doasha}.
\item
If we run the code used to obtain the results displayed in this section and in~\cite{roszas}
several times (having at first settled on whether or not to use TV regularization),
different solutions $\mm$ are obtained because of the different way that the noisy data sets are sampled.
The {\em variance} of such solutions may provide information, or indication, on the extent to which,
albeit having done our best to control the misfit function $\phi (\mm)$, we can trust the quality of
the reconstruction~\cite{sdr}. A large variance in $\mm$ suggests a capricious dependence of the solution on the
noisy data (that is given indirectly through $F(\mm)$), and it would decrease our confidence in
the computed results using a particular bias such as TV regularization.
Having such information is an advantage of the stochastic algorithm over a similar determinstic one.  
\end{itemize}

\section{Conclusions and further notes}
\label{sec:conclusions}

Mathematical software packages typically offer a default option for the error tolerances
used in their implementation. Users often select this default option without much further thinking,
at times almost automatically.
This in itself suggests that practical occasions where the practitioner does not really have
a good hold of a precise
%deterministic 
tolerance value are abundant.
However, since it is often convenient to assume having such a value, 
and convenience may breed complacency, surprises may arise.
We have considered in Section~\ref{sec:examples}
four case studies which highlight various aspects of this uncertainty in 
a tolerance value for a stopping criterion. 

Recognizing that there can often be a significant uncertainty regarding the actual 
tolerance value and the stopping criterion,
we have subsequently considered
%, again using case studies, 
the relaxation of the setting into a probabilistic one, and demonstrated its benefit in
the context of the case study of Section~\ref{sec2.4}. 
The environment defined by~\eqref{3.1} or~\eqref{prob_ineq_lower_upper}, although well-known in other
research areas, is relatively new
(but not entirely untried) in the numerical analysis community. 
It allows, among other benefits, specifying an amount of trust in a given tolerance
using two parameters that can be tuned, as well as the development of bounds on the
sample size of certain Monte Carlo methods, as described in Section~\ref{sec:prob}. 
In Section~\ref{sec:application} we have then
applied this setting in the context of a particular inverse problem involving
the solution of many PDEs, and we have obtained some uncertainty quantification
for a rather efficient algorithm solving a large scale problem.

There are several aspects of our topic that remain untouched in this article.
For instance, there is no discussion of the varying nature of the error quantity that is being
measured (which strongly differs across the subsections of Section~\ref{sec:examples},
from solution error through residual error through data misfit error for an ill-posed problem to 
stochastic quantities that relate even less closely to the solution error). 
Also, we have not mentioned that complex algorithms often involve sub-tasks such as solving a 
linear system of equations iteratively, 
or employing generalized cross validation (GCV) to obtain a tolerance value, 
or invoking some nonlinear optimization routine,
which themselves require some stopping criterion: 
thus, several occurrences of tolerances in one solution algorithm are common.
In the probabilistic sections, we have made the choice of concentrating on bounding the sample size $n$
and not, for example, on minimizing the variance as in~\cite{hutchinson}.

What we have done here is to highlight an often ignored yet rather fundamental issue from different angles. 
Subsequently, we have pointed at and demonstrated a promising approach (or direction of thought) that is not currently 
common in the scientific computing community. Along the way we have also made in context several observations
(Sections~\ref{sec2.1.1},~\ref{sec2.2.1},~\ref{sec2.3.1} and~\ref{sec4.1}) 
from an array of fields of numerical computation, observations which we believe are not common knowledge.

%%%%%%%%%%%%%%%%%%%%%%%%%%%%%%%%%%%%%%%%%%%%%%%%%%%%%%%%%%%%%%%%%%%%%%%%%%%%%
\section*{Acknowledgments}
\label{acknowledgment}
The authors thank Eldad Haber and Arieh Iserles for several fruitful discussions. 

%%%%%%%%%%%%%%%%%%%%%%%%%%%%%%%%%%%%%%%%%%%%%%%%%%%%%%%%%%%%%%%%%%%%%%%%%%%%%

%%%%%%%%%%%%%%%%%%%%%%%%%%%%%%%%%%%%%%%%%%%%%%%%%%%%%%%%%%%%%%%%%%%%%%%%%%%%%

\bibliographystyle{plain}
\bibliography{biblio}

\end{document}